\input amstex
\input amsppt.sty
\NoBlackBoxes
\magnification=1200
\parindent 20 pt
 \vsize=7.50in
\NoBlackBoxes
\define\df{\dsize\frac}
\define\bk{\bigskip}

\define \r{\rho}

\define \tX{\tilde{X}}

\define \CP{\Bbb C\Bbb P}

\define \CPt{\Bbb C\Bbb P^2}

\define \vp{\varphi}

\define \Dl{\Delta}

\define\BC{\Bbb C}

\define \1{^{-1}}
\define \2{^{-2}}

\define \Gal{\operatorname{Gal}}

\define \pf{\demo{Proof}}
\define \fc{\frac}

\define \edm{\enddemo}
\define \ep{\endproclaim}

\topmatter

\heading{\bf Chern Classes of Fibered Products of Surfaces}\rm\endheading

\centerline{Mina Teicher}
\abstract\nofrills{\bf Abstract.}
In this paper we introduce a formula to compute Chern classes of fibered
products of
algebraic surfaces.
For $f: X\to \CPt$ a generic projection of an algebraic surface, we define
$X_k$ for
$k\le n$ $(n=\deg f) $ to be the closure of $k$ products of $X$ over $f$
minus the big
diagonal. For $k=n$ (or $n-1)$,\ $X_k$ is called the full Galois cover of
$f$ w.r.t.
full symmetric group. We give a formula for $c_1^2$ and $c_2$ of $X_k.$
For $k=n$ the formulas were already known.
We apply the formula in two examples where we manage to obtain a surface
with a high
slope of $c_1^2/c_2.$
We pose conjectures concerning the spin structure of fibered
products of Veronese surfaces and their fundamental groups.\endabstract

\thanks  This research was partially supported by the Emmy Noether Research
Institute
 of Bar-Ilan University and the Minerva Foundation of Germany.\endthanks
\endtopmatter
\baselineskip 20pt \subheading{0.\ Introduction}

When regarding an algebraic surface $X$ as a topological 4-manifold, it has
the Chern
classes $c_1^2,$\ $c_2$ as topological invariants.
These Chern classes satisfy:
$$\gather c_1^2,\quad c_2\ge 0\\
5c_1^2\ge c_2-36\\
\text{Signature} = \tau=\fc{1}{3}(c_1^2-2c_2)\endgather$$

The famous Bogomolov-Miyaoka-Yau inequality from 1978 (see \cite{Re},
\cite{Mi}, \cite{Y}) states that the Chern classes of an algebraic surface
also satisfy the inequality $$c_1^2\le 3c_2.$$ It is  known that this
inequality is the best possible since Hirzebruch showed  in 1958 that the
equality is
achieved by complex compact quotients of the unit ball (see \cite{H}).

We want to understand the structure of the moduli space of all surfaces
with given
  $c_1^2,\ c_2$; and, in particular, to populate it with interesting
structures of
surfaces.
As a first step it is necessary to develop techniques to compute Chern
classes of
different  surfaces.

In this paper we compute Chern classes of Galois covers of generic
projections
 of surfaces.
 This was already computed in \cite{MoTe2} for the case of the full Galois
cover, where
the product is taken $n$ times ($n$ is the degree of the projection).
In this paper we deal with products  taken $k$ times, $k < n,$ and we
manage to give an
example of a surface where the slope $(c_1^2 / c_2)$ is very high (up to 2.73).
In subsequent research, using the results of this paper and of our ongoing
research on
this subject, we plan to further study these constructions, to compute these
fundamental groups and to decide when the examples are  spin, of positive
index, etc.
We conjecture that for $X_b$ the Veronese surface of order $b$ greater than
$4,$\ $ X_k$ is spin if $k$ is even or $b=2,3(4).$
We further conjecture that for the Hirzebruch surfaces in general the
fundamental groups
of $X_k$ are finite.

In   \cite{RoTe}, we used similar computations to produce a series of
examples of
surfaces with the same Chern classes and different fundamental groups which
are spin
manifolds where one fundamental group is trivial and the other one has a
finite order
which is increasing  to infinity.
The computations in this paper will lead to more examples of pairs in the
$\tau>0$
area.

We consider in this paper fibered products and  Galois covers of generic
projections
of  algebraic surfaces.
If
$f: X\to
\CPt$ is generic of $\deg n,$ we define the $k$-th Galois cover for $k\le
n$ to be
$\underset { f\qquad f\qquad} \to {\overline{X\times\dots\times X-\Delta}}$
where
$\Dl$ is the big diagonal and the fibered product is taken $k$ times. There
exists a
natural projection $g_k:X_k\to\CPt,$ $\deg g_k=n(n-1)\dots (n-k+1).$

The surface $X_k$
for $k=n$,  is called the full Galois cover  (i.e., the Galois cover w.r.t.
full
symmetric group), and is also denoted    $X_{\Gal}$ or $\tX$.
Clearly,  $\deg(X_{\Gal}\to
\CPt)=n!.$ It can be shown that $X_n\simeq X_{n-1}$.
The full Galois covers were first treated by Miyaoka in \cite{Mi}, who
  noticed that their signature should be positive.
 In our papers \cite{MoTe1}, \cite{MoTe2},
\cite{MoTe3}, \cite{MoRoTe}, \cite{RoTe}, \cite{Te}, \cite{FRoTe}, we
discussed the
full Galois covers for $X=f_{|a\ell_1+b\ell_2|}(\CP^1\times\CP^1),$
Veronese embeddings
and Hirzebruch surfaces. In the  papers cited above we computed their
fundamental
groups (which are finite), the Chern  numbers and the divisibility of the
canonical
divisor (to prove that when considered as 4-manifolds they are spin
manifolds).
$X_{\Gal}$ are minimal smooth surfaces of general type.
Other examples of interest on surfaces in the $\tau>0$ area can be found in
\cite{Ch}
and \cite{PPX}.

\bigskip

\subheading{\S1.\ The Main Theorem}

We start with a precise definition.
\definition{Definition}  {\bf A Galois cover of a generic projection w.r.t. the
symmetric group} $\bold {S_k}$ {\bf(for}
$\bold{ k<}$ {\bf degree of the generic projection).}
 Let $ X\hookrightarrow  \CP^N$  be an embedded algebraic surface.
Let $f: X\to\CPt$  be a generic projection, $n=\deg f.$
 For $1\le  k\le n,$\
 let
 $$\align
&X\underset{f\quad f}\to{\times\dots\times}
X=\{(x_1,\dots x_k)\bigm |x_i\in  X,\ f(x_i)=f(x_j)\ \forall i\forall j\} ,\\
&\Dl=\{(x_1,\dots,x_k)\in  X\underset f\quad
f\to{\times\dots\times} X\bigm| x_i=x_j\ \ \text{for
some}\ i\ne j\}\\
& X_k= \underset{\underset k\to{\underbrace{\ \  f\qquad
f }}}\to{\overline{X \times\dots\times}}
\overline{X-\Dl}.\endalign $$
  $X_k$  is the closure of $
X\underset{ f\quad f}\to{\times\dots\times} X-\Delta$.
 \ $X_k$  is the  Galois cover w.r.t.
the symmetric group on $k$ elements.
 We denote $X_0 = \CPt.$ \enddefinition
\smallskip
For every $k \ge 1$ we have the canonical projections $g_k:X_k\to\CPt$ and
a natural
projection (on the first $k$ factors)
$f_k:X_k\to X_{k-1},$  which satisfy\newline
 $$f_1=g_1=f$$
$$g_{k-1}f_k=g_k, \quad(k \ge 2).$$
Clearly, $$\align &\deg g_k=n\cdot (n-1)\dots (n- k+1)\\ &  \deg f_k=n-k+1,\\
&X_{n-1}\simeq X_n\ \text{($f_n$ is
an isomorphism).}\endalign$$.

For $k=n$ \ (or $n-1)$, we call $X_k$ the Galois cover w.r.t. the full
symmetric group
or the full Galois cover and denote it also by $X_{\Gal}$.

\remark{{\bf Remark}} $X_k$ is the interesting component in the fibered
product
$\underset{\underset k\to{\underbrace{  f\qquad
f }}}\to{X \times\dots\times X}$
\endremark

\demo{Notations}

For the rest of the paper we shall use the following notations:

$n=\deg f.$

$X_k,$ the Galois cover  of $f: X\to\CPt$ as above, $k\le n.$

$S=$ the branch curve of $f$ in $\CPt$ ($S$ is a cuspidal curve)

$m=\deg S$

$\mu=\deg S^* $ $(S^*$ the dual to $S)$

 \quad  = number of branch points in $S$ w.r.t. a generic projection of
$\BC^2$ to
$\BC^1.$

$d =$ number of nodes in $S$

$\rho =$ number of cusps in $S$\enddemo

\proclaim{Theorem 1} The Chern classes of $X_k$ are as follows:

\flushpar\rom{(a)} \qquad\quad
$c_1^2(X_1)=9n+\left(\df{m}{2}-6\right)m-\rho-d.$

 For $2\le k\le n- 1$
$$\align c_1^2(X_k)&=9(n-k+1)\dots n\\
&+\frac{1}{2}[(n-k+1)\dots (n-2)](2n-k-1)k\left(\frac{m}{2}-6\right)m\\
&-[(n-k-1) \dots (n-3)]k\rho\\
&-\fc{1}{2}[(n-k-1)\dots (n-4)](2n-k-5)kd\endalign$$
\flushpar\rom{(b)}
$$\align c_2(X_1)&=3n-2m+\mu\\
c_2(X_2)&=3n(n-1)-2(2n-3)m+(2n-3)\mu+\rho+2d\\
c_2(X_3)&=3n(n-1)(n-2)-3(2n-4)(n-2)m
 +\fc{3}{2}(2n-4)(n-2)\mu\\
&+2(3n-9)d+(3n-8)\rho\endalign$$

For $  4\le k\le n-1$
$$\align c_2(X_k)&=3(n-k+1)\dots n\\
&   -(n-k+1)\dots (n-2) (2n-k-1)k  m\\
&+\frac{1}{2}(n-k+1)\dots (n-2)  (2n-k-1)k\mu\\
&+(n-k+1)\dots(n-3)(k-1)k\left(\frac{n}{2}-\frac{k+1}{3}\right)\rho\\
&+[(n-k+1)\dots(n-4)]\fc{k(k+1)}{4}\{(k+6)(k-1)+4n(n-k-1)\}d\\
&+[(n-k+1)\dots (n-4)]\{4nk-2n^2k\}d\endalign$$
\endproclaim

\remark{{\bf Remarks}}

(a) We consider an empty multiplication as 1.

(b) The case $k=n-1 (X_k=X_{\Gal}),$ of this Theorem was treated in
\cite{MoTe2},
Proposition 0.2 (proof there is given by  F. Catanese).  (See also
\cite{MoRoTe}).
One can easily see that for $k = n$ the formulas here coincide with the
formulas from
\cite{MoTe2}.
For $c^2_1$ it is enough to use remark (a) about empty multiplication.
We get:
$$c_1^2(X_{\Gal})=c_1^2(X_{n-1})=\frac{n!}{4}(m-6)^2 .$$
Note that $d$ and $\rho$ do not appear in this formula.
For $c_2$ we get here (using $(a_1)$)
$$c_2(X_{\Gal})=c_2(X_{n-1})=n!\left(3-m+\frac{1}{4}d+\frac{\mu}{2}+\frac{\rho}{
6}
\right)$$
which coincide with \cite{MoTe2}, using the formula for the degree of
the dual curve:
$$\mu = m^2 - m - 2d - 3\rho.$$\endremark

\demo{Proof of the Theorem}

Let $g_k: X_k\to\CPt,$\ $f_k: X_k\to X_{k-1}$ be the natural projections.
Clearly,  $g_1=f_1=f,$\ $g_k=g_{k-1}f_k$, for $k\ge 2$\quad $\deg
f_k=n-k+1,$\ $\deg
g_k=\dsize\frac{n!}{(n-k)!}.$
Let $E_k$ and $K_{X_k}$ be the  hyperplane and canonical divisors of $X_k,$
respectively ($E_k = g^*(\ell)$ for a line $\ell$ in $\CPt$).

Let $S_k$ be the branch curve of $f_k$ (in $X_{k-1}),$ \ $m_k$ its degree
and $\mu_k'$
the number of branch points that do not come from $S_{k-1}$\ $ (S_1 = S).$
Let $S_k'$
be the ramification locus of
$f_k$  (in
$X_k).$ Let $T_k'$  be the ramification locus of $g_k$  (in $X_k).$

We recall that the branch points in $S$ (or $S_k$) come from two points
coming together
in the fibre, the cusps from (isolated) occurrences of three points coming
together and
nodes from 4 points coming together into 2 distinct points.
Generically, cusps and nodes are unbranched.
We use this observation in the sequel.

To compute $c_1^2(X_k)$ we shall use:
$$\align &c_1^2 (X_k)=K_{X_k}^2\\
&K_{X_k}=-3E_k+T_k'\endalign$$
and the following identities.
$$\align
&T_k'=\left\{\aligned &S_k'+f_k^*(T_{k-1}')\qquad k\ge 2\\
&S_1'\qquad\qquad\qquad\quad k=1\endaligned\right.\\
&T_k'=-\frac{1}{2}S_{k+1}+\frac{1}{2}g_k^*(S)\endalign$$

To compute $c_2(X_k)$ we shall assume that all cusps and nodes of $S$ are
vertices of a
triangulation.
Using the standard stratification computations, this implies the following
recursive formula:
$$c_2(X_k)=\deg f_k\cdot c_2(X_{k-1})-2m_k+\mu_k'.$$
Thus we need to get a formula for $E_k\cdot T_k',$\quad $S_{k+1}\cdot
T_k',$\quad $m_k $
and $\mu_k'.$ We shall use the following 3 claims:\enddemo

\proclaim{Claim 1}

\rom{(i)}\ \ Let $m_k=\deg S_k.$ For $k\ge 2,$\ $m_k =(n-k)\dots
(n-2)m,$\quad $m_1=m.$

\rom{(ii)}\ \  Let $d_k=$ \# nodes in
$S_k.$ For $k\ge 2$, \ $d_k= (n-k-2)\dots (n-4)d,$\quad $d_1=d$.

\rom{(iii)}\  Let $\rho_k=$ \# cusps in $ S_k.$ For $k\ge 2,$\
$\rho_k=(n-k-1)\dots
(n-3)\rho,$\quad $\rho_1=\rho.$

\rom{(iv)}\ Let $\mu_k'$ be the number of branch points of $S_k$ that do
not come from
$S_{k-1},$\linebreak
$\mu_k'=
\mu_k-(n-k+1)\mu_{k-1}$\ $(k\ge 2)$ and   $\mu_1'=\mu.$
Then for $k=2$,
 $\mu_2'=(n-2)\mu+\rho+2d$ and for  $k\ge 3$\ $\mu_k'=(n-k)\dots (n-2)\mu+$
$(n-k)\dots (n-3)(k-1)\rho+[(n-k)\dots $\linebreak$(n-4)](k-1)(2n-k-4)d.$
(For $k=3$ the coefficient of $d$ is $2(2n-7).)$
\endproclaim

\proclaim{Claim 2}
$$E_{k.}T_k'=\cases m&\quad   k=1\\
 \frac{1}{2} m[(n-k+1)\dots (n-2)]
\{(2n-1)k-k^2\}&\quad k\ge 2\endcases$$
\endproclaim

\proclaim{Claim 3}
$$S_{k+1.}T_k'=\cases 2\rho +2d\quad & k=1\\
2(n-k-1)\dots (n-3)k\rho+(n-k-1)\dots (n-4)(2n-k-5)kd\quad & k\ge 2\endcases$$
\endproclaim

\demo{Proof of Claim 1}

Items (i), (ii) and (iii) are easy to verify from the definition of fibered
product.
For (iv) we notice that $\{\mu_k'\}$ satisfy the following recursive equations:
$$\align&\mu_k'=(n-k)\mu_{k-1}'+\rho_{k-1}+2d_{k-1}\quad k\ge 2\\
&\mu_1'=\mu.\endalign$$

The formula for $\mu_2',\mu_3'$ follows immediately from the recursive formula.
For $k\ge 4$ we
substitute the formulas for
$\rho_{k-1}$ and
$d_{k-1}$ from (ii) and (iii) to get
$\mu_k'=(n-k)\mu_{k-1}'+ (n-k)\dots (n-3) \rho+2 (n-k-1)\dots (n-4) d$ and
we shall
proceed by induction.
By the induction hypothesis $\mu_{k-1}'= (n-k+1)\dots (n-2) \mu+ (n-k+1)\dots
(n-3) (k-2)\rho+ [(n-k+1)\dots(n-4)](k-2)(2n-k-3)d.$

\flushpar We substituted the last expression in the previous one to get
$$\align\mu_k'&=(n-k)(n-k+1)\dots (n-2)\mu+(n-k)(n-k+1)\dots(n-3)(k-2)\rho \\
&+(n-k)(n-k+1)\dots
(n-4)(k-2)(2n-k-3)d\\
&+(n-k)\dots(n-3)\rho+2(n-k-1)\dots (n-4)d\\
&=(n-k)\dots(n-2)\mu+(n-k)\dots(n-3)(k-1)\rho\\
&+(n-k)\dots(n-4)\{(k-2)(2n-k-3)+2(n-k-1)\}d\endalign$$
which coincide with the claim since $(k-2)(2n-k-3)+2(n-k-1)=$\linebreak
$(k-1)(2n-k-4).$
\hfill    $\qed$ for Claim 1
\enddemo

\demo{Proof of Claim 2}

For $k\ge 2$
$$\allowdisplaybreaks\align E_{k.} T_k'&=\frac{1}{2} E_{k.}(g_k^*(S)-S_{k+1})\\
&=\frac{1}{2}E_kg_k^*(S)-\frac{1}{2}E_{k.} S_{k+1}\\
&=\frac{1}{2}g_k^*(\ell)  g_k^*(S)-\frac{1}{2}E_{k.} S_{k+1}\\
&=\frac{1}{2}g_k^*(\ell. S)-\frac{1}{2}E_{k.}S_{k+1}\\
&=\frac{1}{2}(\deg g_{k.}) m-\frac{1}{2}m_{k+1}\\
&=\frac{1}{2}m(n-k+1)\dots n-\frac{1}{2}(n-k-1)(n-k)\dots(n-2)m\\
&=\frac{1}{2}m[(n-k+1)\dots (n-2)]\{(n-1)n-(n-k-1)(n-k)\}\\
&=\frac{1}{2}m[(n-k+1)\dots (n-2)]\{2nk-k-k^2\}.\qquad\qquad\qquad \qed\
\text{for
Claim 2}\endalign$$

\enddemo

\demo{Proof of Claim 3}

Since $T_1'=S_1',$ the formula trivializes for $k=1.$
$S_2\cdot T_1'=S_2\cdot S_1'=2\rho+2d=2\rho_1+2d_1.$
For $k\ge 2$
$$\align
S_{k+1}\cdot T_k'&=S_{k+1}(f_k(T'_{k-1})'+S_k')\\
&=S_{k+1}\cdot f_k^*(T_{k-1}')+(S_{k+1}\cdot S_k')\\
&=(\deg f_k\bigm|_{S_{k+1}})\cdot (S_k\cdot T_{k-1}')+2\rho_k+2d_k\\
&=(\deg f_k -2)(S_k\cdot T_{k-1}')+2\rho_k+2d_k\\
&=(n-k-1)(S_k\cdot T_{k-1}')+2\rho_k+2d_k. \endalign $$
Denote $a_k=S_{k+1}\cdot T_k'$.

We shall prove the claim by induction using the recursive formula
\linebreak $a_k=
(n-k-1)a_{k-1}+2\rho_k+2d_k$.
  For $k=2:$
$$\align a_2&=(n-3)a_1+2\rho_2+2d_2\\
&=(n-3)(2\rho+2d)+2(n-3)\rho+2(n-4)d\\
&=4(n-3)\rho+2d(n-3+n-4)\\
&=4(n-3)\rho+2d(2n-7).\endalign$$
Thus the statement is true for $k=2.$

Let $k\ge 3.$ Assume the formula is true for $k-1.$ We shall prove it for $k.$
$$\align
a_k&=(n-k-1)a_{k-1}+2\rho_k+2d_k\\
&=(n-k-1)\{2(n-k)\dots(n-3)(k-1)\rho+(k-1)(n-k)\dots(n-4)(2n-k-4)d\}\\
&+2(n-k-1)\dots(n-3)\rho+2(n-k-2)\dots(n-4)d\\
&=2(n-k-1)\dots(n-3)k\rho+(n-k-1)\dots(n-4)\{(2n-k-4)(k-1)+2(n-k-2)\}\\
&=2(n-k-1)\dots(n-3)k\rho+(n-k-1)\dots(n-4)\{2nk-k^2-5k\}d.
\endalign$$
>From the two formulae, we can see that the product $(n-k-1) ... (n - 4)$
should be
1 for
$k\le 2.$
\hfill    $\qed$ for Claim 3
\enddemo

We go back to the proof of the theorem.
To prove (a) we write
$$\align c_1^2(X_k)&=K_{X_k}^2=(-3E_k+T_{k}')^2\\
&=9E_k^2-6E_k. T_k'+(T_k')^2\\
&=9E_k^2-6E_{k.}T_k'+T_k'\left[-\frac{1}{2}S_{k+1}+\frac{1}{2}g_k^*(S)\right]\\
&=9E_k^2-6E_k\cdot T_k'-\frac{1}{2}T_k'\cdot S_{k+1}+\frac{1}{2}T_k'\cdot
g_k^*(S).\endalign$$

Now: $E_k^2=\deg g_k=(n-k+1)\dots n$ $\left(=\dsize\frac{n!}{(n-k)!}\right).$
Since $S$ is of $\deg m$\ $T_k'. g_k^*(S)=mE_k. T_k'.$
We substitute the results from Claim 2 and Claim 3 to get (a).

We prove (b) by induction on $k.$
For $k=1$ we take the recursive formula
$c_2(X_k)=(n-k+1)c_2(X_{k-1})-2m_k+\mu_k'$ and
substitute $k=1$ to get $c_2(X_1) =3n-2m+\mu$ which coincides with formula
(b) for
$k=1.$ We do the same for $k=2,3.$
To prove $k-1$ implies $k$ we use Claim 1(iv) and (ii) to write
$$\align c_2(X_k)=(n-k+1)&c_2 (X_{k-1})-2(n-k)\dots (n-2)m+(n-k)\dots
(n-2)\mu\\
&+(k-1)(n-k)\dots (n-3)\rho\\
&+(n-k)\dots (n-4)(k-1)(2n-k-4)d\endalign$$
When substituting the inductive statement for $c_2(X_{k-1})$ and shifting
around
terms,  we get (b).

\hfill    $\qed$ for the Theorem

\bigskip
\subheading{\S2.\ A Different Presentation of the Chern Classes}
\proclaim{Proposition 2}
 Let $E$ and
$K$ denote    the hyperplane and canonical divisors of
$X$, respectively.  Then the Chern classes of $  X_k$ are functions of
$c_1^2(X)$,\ $c_2(X)$,\ $\deg (X)$,\   $E,$\ $ K,$ and $k$.\endproclaim

\demo{Proof} (Proof for $X_n$ appeared in \cite{RoTe})
Let $S$ be the branch curve of the generic projection $f: X\to\CPt$\
$(S\subseteq
\CPt)$.
By Theorem 1, the Chern classes of $X_k$ depend on $k,$\ $\deg(S),$
$\deg(X)$ and
$\mu,d,\rho,$ the number of branch points, nodes and cusps of $S,$
respectively.

We shall first show   that $\mu,d,\rho$ depends on $c_2(X),$\  $\deg X,$\
$\deg(S),$\
$e(E)$ and
$g(R)$ where $g$ denotes the genus of an
algebraic curve,   $e$ denotes the topological Euler characteristic
of a space, and $R$\ $( \subset X)$ is the ramification locus of $f$ which
is, in
fact, the  non-singular model of $S$.

Recall that $\mu$ also is equal to $\deg(S^*),$ where $S^*$ is the dual
curve to $S.$
For short we write $n=\deg(X),$\ $m=\deg(S).$

  We show this by presenting three linearly independent formulae:
$$\align &\mu=m(m-1) - 2d - 3 \rho\\
&g(R) = \frac {(m-1)(m-2)} 2 - d - \rho\\
&c_2(X) + n = 2e(E) + \mu\endalign$$
The first two are well-known formulae for the degree of the dual curve and
the genus
of a non singular model of a curve.  For the third, we may find a Lefschetz
pencil of
hyperplane sections of
$X$ whose union is
$X$.
The number of singular
curves in the pencil is equal to $\mu$.
The topological Euler characteristic of the fibration equals  $e(X)  = e(\CP^1)
\cdot e(E) +
\mu-n $ ($n$ appears from blowing up $n$ points in the hyperplane sections).
The formula follows from $e( \CP^1)=2$ and $e(X)=c_2(X).$

We shall conclude by showing  that $\deg(S)$,\ $e(E)$ and $g(R)$ depend on
$c_1^2(X),$\
$\deg X$ and $E.K.$

This follows from the Riemann-Hurwitz formula, $R=K+3E,$ the adjunction formula
$2-2g(C)=-C.(C+K),$ and the fact that $E^2=\deg X$ and $K^2=c_1^2(X).$
In fact, we have:
$$\align &g(R)=1+\fc{1}{2}R(R+K)=1+\fc{1}{2}(K+3E)(2K+3E)\\
&e(E)=2-2g(E)=-E(E+K)\\
&\deg(S)=\deg(R)=E.R=E(K+3E). \qquad\qquad\qquad\qquad\qed\endalign$$
 \enddemo

>From the above proof we can, in fact, get the precise formulae of
$c_1^2(X_k)$ and
$c_2(X_k)$ in terms of $c_1^2(X),$\ $c_2(X),$\ $\deg(X)$, $E.K,$ and $k.$
For certain (computerized) computations,  it is easier to work with these
formulae
rather than those of \linebreak Theorem 1.

\proclaim{ Corollary 2.1}
In the notations of the above  proposition:
$$\align  &c_1^2(X_{n})=\frac{n!}{4}[(E. K)^2+6n(E. K)+9n^2-12(E.
K)-36n+36]
\\  &c_2 (X_{n})=\\ &\frac{n!}{24}[72-10c_1^2(X)-54(E.
K)-114n+27n^2+14c_2(X)+3(E. K)^2+18n(E. K)]\endalign$$\endproclaim

Similar formulas can be obtained for $X_k$ for $k<n.$
\bigskip

\bk
\subheading{\S3.\ Examples}

To use Theorem 1, we need computations of $n,$ $m$, $\mu,$ $\r$ and $d.$
We compute them for two examples.

\proclaim{Examples 3.1}
For $X=V_b$, a Veronese embedding of order b, we have
$$\align & n=b^2\\
&m=3b(b-1)\\
&\mu=3(b-1)^2\\
&\vp=3(b-1)(4b-5)\\
&d=\df{3}{2}\ (b-1) (3b^3-3b^2-14b+16)\endalign$$
(see \cite{MoTe3}).\ep

\demo{Proof}\
For $ n, m, \mu$ and $\r,$ see \cite{MoTe3} and \cite{MoTe4}.
 Since $\mu=m(m-1)-2d-3\rho,$ we get the following formula for $d\:$
$2d=m^2-m-\mu-3\vp$ and thus
$$\align 2d&=3b(b-1)(3b(b-1)-1)-3(b-1)^2-9(b-1)(4b-5)\\
&=3(b-1)\{(3b^2-3b-1)b-(b-1)-3(4b-5)\}\\
&=3(b-1)\{3b^3-3b^2-14b+16\}.\endalign$$
\enddemo

When one substitutes $b=3$ and $k=4$, one gets $\fc{c_1^2}{c_2}=2.73.$
By experimental substitutions it seems that  for large b, the signature
$\tau(X_k)$
\ ($= c_1^2-2c_2$), changes from negative to positive at about $\fc{3}{4}n.$

\proclaim{Example 3.2} For
$X=X_{t(a,b)}=f_{|a\ell+bC_+|}$ (Hirzebruch surface of order $t),$
	where  $\ell$ is
a fiber,   $(C_+)^2=t,$  and
$a\ge 1$, we have
$$\align &n = 2ab+tb^2\\
&m=6ab-2a-2b+t(3b^2-b)\\
&\mu=6ab-4a-4b+4+t(3b^2-2b)\\
&\vp=24ab-18a-18b+12+t(12b^2-9b)\endalign$$\ep

\demo{Proof} \cite{MoRoTe},\ Lemma 7.1.3.\edm
\proclaim{Example 3.3} (in the $\tau<0$ area)

\flushpar For $X $ a $K3 $ surface:

$K=0$

$c_1^2(X)=K^2=0$

$c_2(X)=24$

$n=4$

$m=12$

$\mu=36$

$\rho=24$

$d=12$

$c_1^2(X_2)=48$

$c_2(X_2)=144$

$c_1^2(X_3)=c_1^2(X_{\Gal})=216$

$c_2(X_3)=c_2(X_{\Gal})=240.$\ep

\pf
It is well known that for a $K3$ surfaces $K=0,$\ $c_1^2=0,$ \ $c_2=24,$\
$S'=3E,$
\ $n=E^2=4.$
Using this we can get $m$ and $\mu$:
$$\align m&=S'.E=3E.E=3E^2=12\\
 \mu&=c_2(X)-2e(E)+n\quad\text{(see proof of Proposition 2)}\\
&=c_2(X)-2(2-2g(E))+n\\
&=c_2(X)+2E.(E+K)+n=36\endalign$$
Now from
$$m(m-1)=\mu+3\rho+2d$$ and
$$\align m(m-3)&=2g(S')-2+2\rho+2d=(K+S'),S'+2\rho+2d\\
&=3E.3E+2\rho+2d=9E^2+2\rho+2d,\endalign$$
we get $2m=\mu-9E^2+\rho=24$ and $\rho=2m+9E^2-\mu=24.$

Moreover, we get $d=\fc{1}{2}(m(m-1)-\mu-3\rho)=12.$
We substitute these quantities in the formula from Theorem 1 to get the
values of the
Chern classes.
\edm

\remark{\bf{Remark}}
For $t=0,$\ $X_{t(a,b)}$ are actually embeddings of $\CP^1\times\CP^1$.
In  \cite{FRoTe}, we computed the fundamental group of $X_n=X_{\Gal}$ for
$X=X_{t(a,b)}$
which is $\Bbb Z_c^{n-2}$ for $c=g.c.d.(a,b).$ Thus for $(a,b)=1$ these
surfaces are
simply connected.
All these surfaces are smooth minimal surfaces of general type.
For $a\ge 6$,\ $b\ge 5,$ the signature of these surfaces is positive.
For five pairs of $(a,b)$, these surfaces have signature 0 (see \cite{MoRoTe}).
Four of these surfaces are simply connected and the fifth one for which
$a=b=5,$
$\pi_1(X_{\Gal})=\Bbb Z_5^{48}$.

In our ongoing research, we shall apply Theorem 1 and Proposition 2   in
order to obtain more examples of non diffeomorphic surfaces or surfaces in
different
deformation families  with the same
$c_1^2 $ and $c_2$, as well as to compute the slope $\fc{c_1^2}{c_2}$ and
to search
for higher slopes.\endremark

We are also interested in the fundamental groups (in particular, in the
finite ones)
and the divisibility of the canonical
class (in particular, the case where the canonical class is divided by 2,
i.e., the
spin case), which we will investigate in a subsequent paper.
The results in this paper are a basis for producing interesting examples of
surfaces
with positive index, $(c_1^2-c_2)$, finite fundamental groups and spin ($K$
even)
structure. In particular, we plan to prove the following two conjectures.
\proclaim{Conjecture}
For $X=V_b$, Veronese of order $b,$\quad $b>4,$ we have  $X_k$ is a spin
manifold
$\Leftrightarrow k$ even or $b=2,3(4).$\ep
\proclaim{Conjecture} For $X=F_{t(a,b)}$ (the Hirzebruch surface),
$\pi_1(X_k)$ is
finite.\ep

\end